\documentclass[11pt]{amsart}
\usepackage[text={6.1in,8.5in},centering]{geometry}
\usepackage[utf8]{inputenc}

\usepackage{amsfonts,amssymb,amscd,graphicx,amsmath,amsthm}
\usepackage{graphicx}

\usepackage[usenames, dvipsnames]{color}
\title{Parsimonious cones}
\author{Aleksandr Berdnikov}
\newtheorem{theorem}{Theorem}[section]
\newtheorem{cor}[theorem]{Corollary}

\newtheorem{df}[theorem]{Definition}

\DeclareMathOperator\vol{vol}
\DeclareMathOperator\Mix{Mix}

\DeclareMathOperator\Bl{Bl}

\newcommand{\C}{\mathcal{C}}
\newcommand{\B}{\mathcal{B}}

\begin{document}

\maketitle

\begin{abstract}
We construct embeddings of simplicial complexes into a (surface of a) simplicial ball whose triangulation has bounded degrees and low volume. This construction can be used either to efficiently ``simplify a complicated space'' by realizing it as a part of a ball/sphere, or to ``complexify'' a sphere --- to give it a specific metric that inherits desired properties from an embedded complex. 
\end{abstract}

\section*{Introduction}

The construction this article focuses on is more of a trick that can be realized and optimized in a multitude of ways depending on the exact situation, giving respectively differing results. However, the core idea can be illustrated in the proof of the following result.

\begin{theorem}
\label{celltree}
Let $X^n$ be a simplicial complex with $N$ simplices. Then for any $m\geqslant 2n+2$ there is a triple
$(\mathcal{B}^m,\mathcal{C}^{n+1},\mathcal{X}^n)$ of simplicial complexes with bounded degrees so that 
\begin{enumerate}
\item $\mathcal{B}\cong B^m$, $m$-ball,
\item $(\mathcal{C},\mathcal{X})$ is homtopy equivalent to a pair of $X$ and a cone $C(X)$ over $X$,
\item $\mathcal{X}\subset \partial \mathcal{B} \cong S^{m-1}$,
\item number of simplices $|\mathcal{C}|$ and $|\mathcal{B}|$ are no more than $c(m,n) N\log(N)^n$.
\end{enumerate}
\end{theorem}

The idea here is to build a (homotopy) cone $\C$ over $X$ with bounded degrees, and then thicken it to get the ball $\B$. The cone $\C$ can be thought of as a fibration of binary trees over an abstract simplex $\triangle^n$. The leaves of a tree at a given point $p\in\triangle$ are the corresponding points $p_i$ in the simplices $\triangle^n_i$ of $X$. The difficulty is in reshuffling the branches of the tree as $p$ moves across $\triangle$, so that when $p$ approaches a face, the right $p_i$'s get grouped together (those, that merge in $X$ at that face). 

To make this ``fibration'' possible we need to align all the simplices of $X$ so that they can be thought as covering the faces of the abstract simplex $\triangle$. That is, we assume that the vertices of $X$ are colored in $n+1$ colors (so that no simplex different vertices of the same color). That can be ensured by passing to the barycentric subdivision $X_{bar}$ of $X$ (a vertex $v\in X_{bar}$ is colored according to the rank of the face with center $v$). We will also assume that every cell of $X$ is on the boundary of some top-dimensional cell. For example, we can attach $n$-dimensional cell to every cell of $X$ that is not.

%{\bf Remark:} the ``vertex'' of the cone $\mathcal{C}$ is represented by an $n$-simplex with side length $\sim \log (N)$ and it is located in $\partial \mathcal{B} \cong S^{m-1}$ opposite $\mathcal{X}$.

\section{A construction of a single tree}

\subsection{Motivation}\
\label{motive}

The theorem~\ref{celltree} for $n=0$ is easily solved by the following example that motivates the general construction. Given a collection of points $X^0$ we can build a parsimonious cone $\mathcal{C}$ over $X$ as a binary tree with $\mathcal{X}=X$ as its leaves, and let the ball $\mathcal{B}$ be a tubular neighborhood of $\mathcal{C}$, as illustrated below.

\medskip

\includegraphics[scale=0.5]{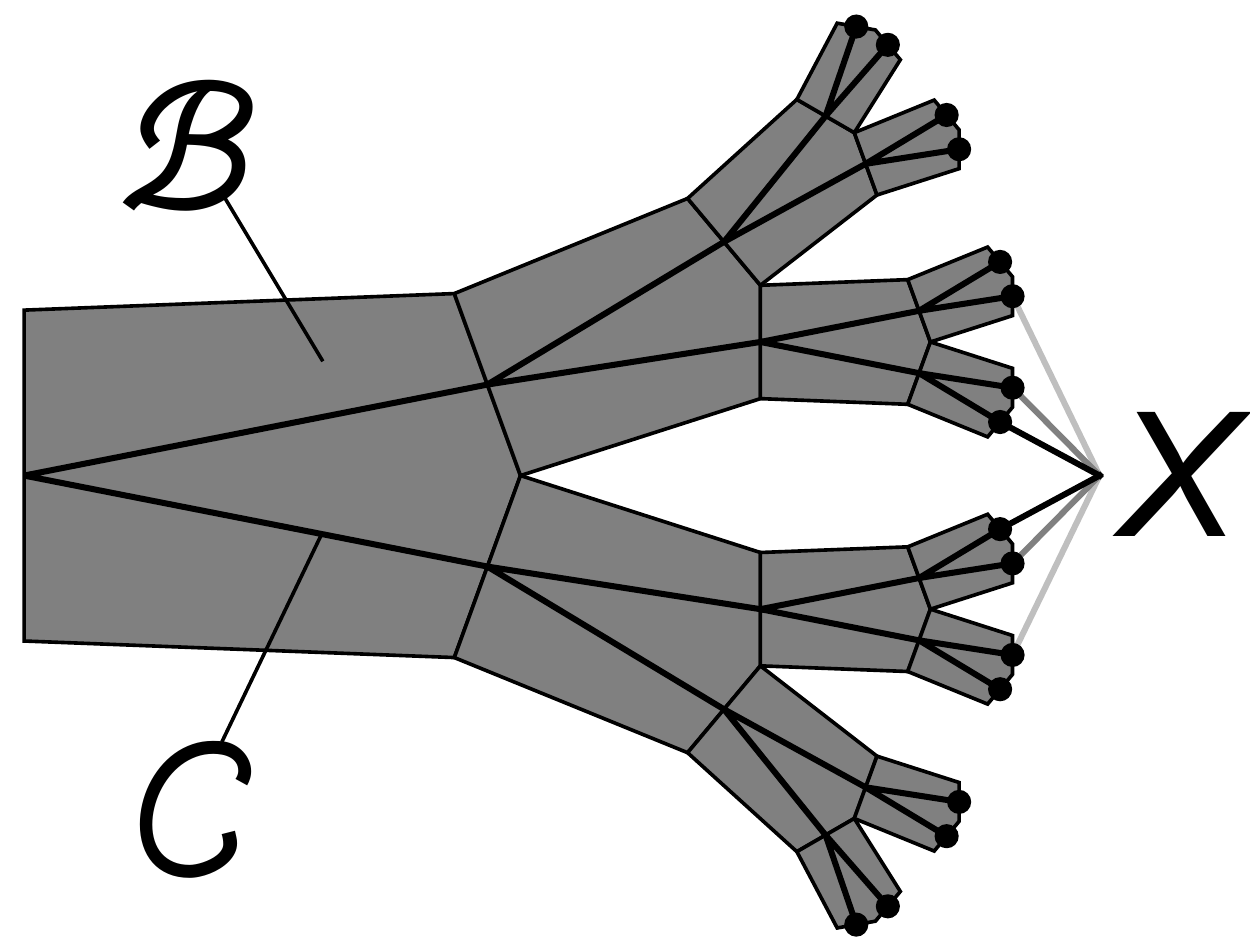}

\medskip

\subsection{Numerations and their mixing}\
\label{numchap}

The structure of the tree can be encoded in enumerating its leaves, each by a string of 0's and 1's describing its location on the tree (say, 0 --- left turn on the split, 1 --- right turn). 

Moving to $n$-dimensional case, we view it, roughly speaking, as an $n$-parametric version of 0-dimensional case. The corresponding parameter space will allow to transition between different tree structures and different leaf numerations.
To work with ``mixed'' numerations that are transitional between the regular ones, we introduce some notation. Assume we have a set $I$ of different numerations $N_{i\in I}:X\rightarrow \{0,1\}^k$ of a (finite) set $X$ by strings of length $k$. Mix them to get numerations, intermediate between $N_i$:

\begin{df}
For strings $S_i\in \{0,1\}^k$ and ``priorities" $p_i\in\mathbb{R}$ assigned to them, a \underline{{\upshape mixed string}} $S_*=\Mix((S,p)_{i\in I})$ is a string of 0's and 1's, each indexed by $i\in I$, written in a line, each with unit spacings, and $S_i$ shifted overall by $p_i$.\\ The \underline{{\upshape index}} of a digit $x$ in $S_*$ is the index $i$ of the string $S_i$ that the digit $x$ came from.
\end{df} 

That is, to get the mixed string, $j$-th digit of $S_i(x)$ is placed at location $j-p_i$ (being indexed by $i$, marking which string this digit came from). An example of string mixing is shown below:

\vskip 0.4cm

\includegraphics[scale=0.3]{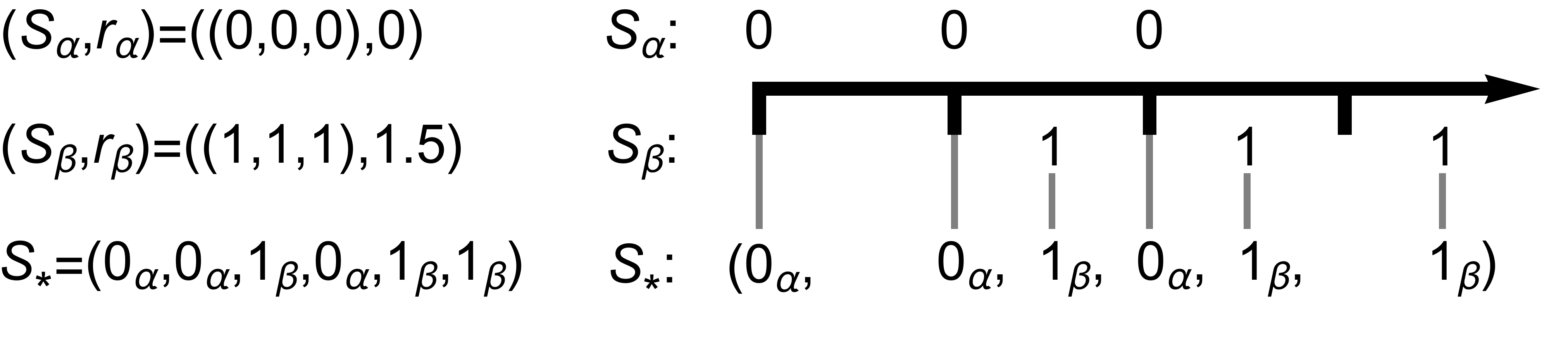}

\vskip 0.6cm

We can mix whole numerations in the same way:

\begin{df}
For numerations $N_i:X\rightarrow \{0,1\}^k$ a \underline{\upshape mixed numeration} $N_*=\Mix((N,p)_{i\in I})$ is defined as mixing of their outputs:
$$N_*(x)=\Mix((N,p)_{i\in I})(x):=\Mix((N(x),p)_{i\in I}).$$ 
\end{df}

\subsection{Building a tree}\

Now, assume we have a mixed numeration of leaves of the tree to be built. The tree is then constructed in a following way: it has root as a unit $n+1$ cube from which it grows upwards and whenever there is a digit indexed $k$ in position $h$ in the numeration, all the branches that carry both leaves with $0_k$ and leaves with $1_k$ in that place, split in half at height $h$ in $k$-th direction. If there are several digits at the same position --- the tree splits in all those directions at the same time. Such a tree may have long but topologically trivial pieces, if in several junctions the branch only carries leaves on one side and hence doesn't split there. To remedy that, all such branches are shortened to length 1.

\section{Proof of the main theorem}

As shown in the motivational section~\ref{motive}, the construction is rather trivial for 0-complexes $X^0$. However, the case $n=1$ that deals with graphs $X^1$ already exposes the crux of the algorithm while being still easy to visualize. For that reason (and because the case of graphs $X^1$ already has non-trivial applications) we first go through the proof in that case and then discuss how to generalize it to arbitrary $n$.

\subsection{Illustrative case of n=1}

{\bf Proof.} Consider a bipartite graph $X^1$ with sets of vertices $X^{(0)}_1$ and $X^{(0)}_2$. Each of these sets provides a numeration $N_i$ of the edges $X^{(1)}$ of the graph: each edge gets a number that is at its endpoint. We build trees using these numerations, mixed in all possible ways: as one numeration moves past the other, the tree changes accordingly. The (centers of the) leaves of the morphing tree will trace lines that we interpret as the edges of $\mathcal{X}$. The endpoint trees (that is, built for $Mix((N_1,0),(N_2,\log(|X|)))$ and for $Mix((N_1,\log(|X|)),(N_2,0))$) each consists of a big lower tree that splits in on direction and from their leaves stem smaller trees that split in the second direction. By construction, each upper tree gathers exactly those edges at its leaves, that meet in one of vertices of $X$. Thus we can join the edges of $\mathcal{X}$ by lines that follow those upper trees on their side. That concludes the construction of $\mathcal{B}$ and $\mathcal{X}$. The cone $\mathcal{C}$ is just the center core of $\mathcal{B}$ meeting the edges of $\mathcal{X}$ plus additional part joining it to the rest of $\mathcal{X}$ (to its ``vertex trees'') in an obvious way.

Here are the specifics. We may assume $|X^{(0)}_1|=|X^{(0)}_2|$ for convenience. Let the numeration $N(t)$ of the edges $X^{(1)}$ at time $t$ be $$N(t)=\Mix\bigl( (X^{(0)}_1,t),(X^{(0)}_2,|X^{(0)}_i|-t)\bigr),$$ so that it passes from $X^{(0)}_1$ followed by $X^{(0)}_2$ towards being the other way around. Let $\{t_i\}$ be the set of times at which the  difference of the priorities is a half-integer, so that $N(t)$ is properly defined. For each $t_i$ there is a tree $T(N(t_i))=:T(i)$, so that successive trees differ only slightly, and so we are left to morph each of them into the next one.

Consider two subsequent trees $T(i)$ and $T(i+1)$. Their numerations differ by swapping some of neighboring values as they move past each other, which corresponds to the following tree changes:

\includegraphics[scale=0.4]{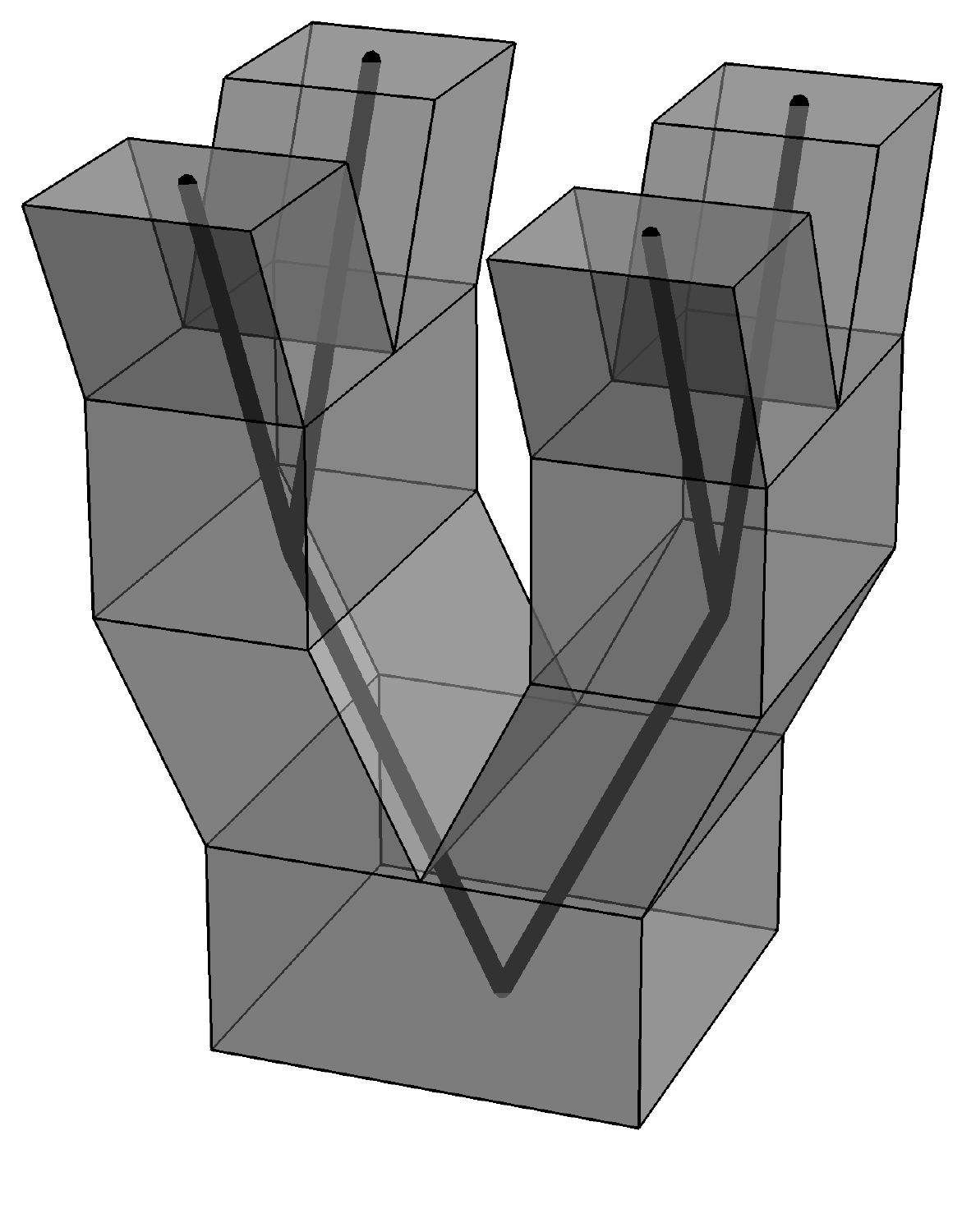} 
\includegraphics[scale=0.4]{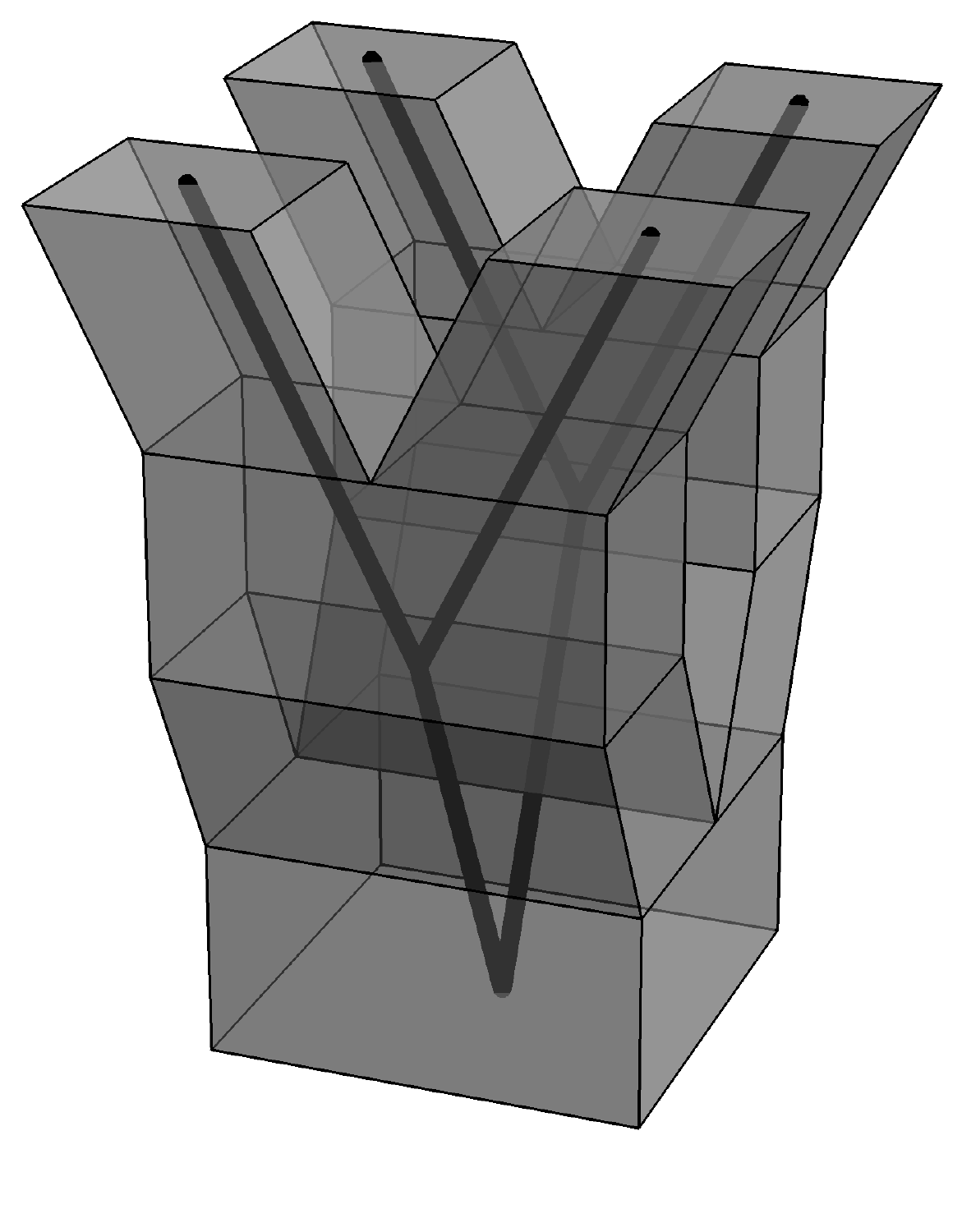}

We exploit their similarity to connect them:

\begin{df}
A \underline{\upshape morphing} between the \underline{\upshape little trees} 
$$T(\cup_{i,j\in\{0,1\}}\overline{ i_\alpha j_\beta})=:T_{\alpha\beta}$$ and $T_{\beta\alpha}$ is a triangulated 4-ball $T_{\alpha-\beta}\cong B^3\times[0,1]$ such that 
\begin{enumerate}
\item the sides $B^3\times\{0,1\}$ of the ball are $T_{\alpha\beta}$ and $T_{\beta\alpha}$; \\
\item the horizontal faces $F_i$ of $B^3\times\{0\}$ are connected by $$F_i\times[0,1]\hookrightarrow \partial B^3\times[0,1]$$ to the analogous faces of $B^3\times\{1\}$;\\
\item the core trees $C_{\alpha\beta}$ and $C_{\beta\alpha}$ are in turn connected by a 2-complex $C_{\alpha-\beta}$ such that $C_{\alpha-\beta}\cap B^3\times\{t\}$ is a tree inside $B^3$ with root and leaves being the middles of the horizontal faces $F_i$.  
\end{enumerate}
\end{df}

If a little tree lacks some parts but still has different splittings that have to be exchanged, the morphing is defined by first growing the little trees back to their full form and morphing them. Thus we can morph any tees by combining the little ones.

\begin{df}
A \underline{\upshape morphing} between the whole trees 
$$T\Bigl(\Mix\bigl((N_1,p_1),(N_2,p_2)\bigr)=:N\Bigr)\text{ and }T\Bigl(\Mix\bigl((N_1,p_1'),(N_2,p_2')\bigl)=:N'\Bigr)$$
(with $(p_1-p_2)-(p_1'-p_2')=1$) is a gluing by the common $F_i\times [0,1]$ faces of morphings between the little trees, that $T(N)$ and $T(N')$ consist of.
\end{df}

Notice that regardless the manipulations, the homotopy type of morphing object is always the same: a triple of a 3-ball; a tree inside, with a root and leaves on the boundary; and those leaves.

Concatentating the morphings --- one for each time the mixed numeration changes --- gives a morphing (of logarithmic length) between the tree that has the leaves grouped according $N_1$ first (so that the edges meet as they do at the first vertices) and a tree grouping by $N_2$ first (so that the edges meet as at the second vertices). Therefore we get the triple $(\mathcal{B},\mathcal{C},\mathcal{X})$ as in the theorem, if we attach cones following the upper trees by their boundary, and connect them to the core $\mathcal{C}$ inside.   

\begin{flushright}
$\square$
\end{flushright}

\subsection{Generalization to arbitrary dimension}

{\bf Proof.} When the complex $X$ is of higher dimension, we just need to interpolate between more numerations corresponding to $n+1$ types of the vertices of the complex. Let the \underline{\upshape base complex} $B$ for $n+1$ numerations be the hyperplane $\mathbb{R}^n\subset \mathbb{R}^{n+1}$ defined by $\sum x_i=0$, with a structure of a complex dual to its Voronoi tesseletion done for the set of points whose coordinates are equidistributed modulo $\mathbb{Z}$.

The base complex parametrizes possible relative priorities of $n+1$ numeration. Its vertices correspond to unambiguous (distinct mod $\mathbb{Z}$) priorities, and higher-dimensional cells are centered at the points where some priorities collide mod $\mathbb{Z}$: these cells connect the points that have these priorities separated in all possible ways. 

The parsimonious cone triple $(\mathcal{B},\mathcal{C},\mathcal{X})$ is constructed mostly as a (homotopy) fibration over the relevant part of $B$, inducting by the skeleta of $B$. The fibers are the binary trees $T(N_i,x_i)$ whose leaves correspond to $n$-dimensional cells of $X$, therefore their size is $\sim N=|X|$. The linear size of the part of $B$ needed (allowing all possible mixing of numerations) is of order of the length of the numerations, that is, $\sim \log (N)$. Thus the whole construction is of size $\sim N \log(N)^n$. 

Over the vertices $B^{(0)}$ with coordinates $(x_i)$ we plant the trees $T(N_i,x_i)$, obtained by mixing the numerations $N_i$ with priorities $x_i$. Having built the construction over the skeleton $B^{(k)}$, we extend it to the cells of $B^{k+1}$ in a following way. By the induction hypothesis we have a tree fibration over the boundary of the cell. Moreover, for any section of the coordinates $x_i$, that is, splitting them into $x_i^{(s)}$ and $x_i^{(g)}$ so that $x_i^{(s)}<x_i^{(g)}$, that splitting marks a floor in the tree that doesn't undergo any changes over this cell, so $\Bl$-faces at that floor extend just by the direct product over the boundary of the cell; so we extend them by the direct product over the cell itself as well. Between those $\Bl$-faces are ball fibrations with a tree inside with the root on the bottom $\Bl$-face and leaves at the top $\Bl$-faces. Homotopically it is a trivial fibration, so we extend it homotopically trivially into the cell. In case at some places of the boundary this fibration underwent shortening, we first attach buffers that lengthen it back. However, if it was shortened over the whole boundary, there in no need for lengthening and the extension is done in a shortened form. That ensures that the whole construction uses the number of simplices on the same order as the starting collection of trees over $B^{(0)}$. Note that the spacing between these splitting $\Bl$-faces can be at most $n+1$ high (that is, how many priorities are clumped up mod $\mathbb{Z}$ and cannot be splitted apart at this cell). Therefore such extensions may be built upon only finite number of different triangulations and topologies on the boundary, so we can again use only finite amount of triangulated extensions to built the fibration over the next dimension.

Now the whole space built is the ball $\mathcal{B}$, the trees inside it form the ``cone'' $\mathcal{C}$, but the demarkation of the $\mathcal{X}$ requires some work. For every $n$-cell of $X$ we include in $\mathcal{X}$ the union of corresponding leaves of the fibers-trees. In total we get disjoint contractible sets corresponding to $n$-cells. Inductively, for $(n-k)$-cell $A$ (with $k>0$), a cone needs to be attached to the link of $A$. Luckily, there is cone attached to everything --- the core tree fibration, --- so we just have to select the right part. The right part is formed by (parts of) core trees fibered over the set of points $(x_i)$, that put the $k$ numerations that aren't present among vertices of $A$, completely on top of all other numerations. The specific part of the fiber is the core subtree that stems from the floor that splits the numerations. The $\Bl$-faces on this floor are differentiated by lower but not upper numerations, so there is a $\Bl$-cell $\Bl(A)$ that carries precisely the leaves to which this cell should be attached. Thus, to the cell $A$ correspond the set $\mathcal{A}$ swiped out by the core trees stemming from $\mathcal{A}$. The contractible fibers of $\mathcal{A}$ and the base that mirrors the cell structure of the simplex secures that homotopically attaching $\mathcal{A}$ to the built so far $\mathcal{X}$ is identical to attaching $A$ by the link of $A$ to the corresponding part of $X$.

\begin{flushright}
$\square$
\end{flushright}

\section{Local simplification of complexes}

One of the things happening in theorem~\ref{celltree} is that the complex $X$ with arbitrary degrees is replaced with a homotopy equivalent complex $\mathcal{X}$ with bounded degrees. If that is the only goal, we can do a little better by one log factor:

\begin{theorem}
Let $X^n$ be a simplicial complex with $N$ simplices. Then there exists a homotopy equivalent complex $\mathcal{X}$ with bounded degrees with $|\mathcal{X}|\sim N \log(N)^{n-1}$ 
\end{theorem}

{\bf Proof} follows from the observation that most of the volume of $\mathcal{X}$ in theorem~\ref{celltree} is carried by leafs fibering over the interior of the base complex $B$. Therefore deflating each such cell so that they no longer have much of interior, the size of $\mathcal{X}$ is reduced to what was fibered over the boundary of the base complex and becomes of size $\sim N \log(N)^{n-1}$.

\begin{flushright}
$\square$
\end{flushright}

There might be a situation where more refined local information about $X$ is known, like the estimate for sizes of links of its cells. In such situation it might be more efficient to build $\mathcal{X}$ not as part of a cone over the whole $X$, but simplifying the links locally. More specifically, $\mathcal{X}$ can be built inducting by the dimension of cells of $X$, starting with the highest. To add a cell, one considers its link --- a complex to which it should be attached --- and builds a parsimonious cone over it using theorem~\ref{celltree}. This way the construction is produced locally and therefore may lead to better results in some situations.

\section{Sphere hard to cut}
\subsection{Setting and preliminaries}

\begin{df}
A space $X$ is said to be $(M,\varepsilon)$-\underline{\upshape hard to cut} if one cannot cut at least $\varepsilon$ volume of $X$ by a cut of area $M$ or less. That is, for any $\Omega\subset X$ with $\vol(\Omega),\vol(X \setminus \Omega)>\varepsilon$, one has $\vol(\partial \Omega)>M$.
\end{df}

The paper~\cite{sphere} presents (for any $M,\varepsilon>0$ and $n\geqslant 3$) a construction of a metric sphere $S^n$ that is $(M,\varepsilon)$-hart to cut. Our technique allows to adjust this construction to produce qualitatively efficient examples. To make the statement precise we introduce a common measure of complexity:

\begin{df}
The complexity of a manifold $(X,g)$ with a metric $g$ on it is $\vol(X,kg)$ where $k$ is the smallest number so that $(X,kg)$ is close to euclidean space on scale 1.
\end{df}

For example if $X$ is closed then $k\sim \max(r^{-1},\sqrt{||R||_\infty})$ where $r$ is injectivity radius of $X$ and $R$ is its Riemann tensor. Such notion of complexity transfers well to simplicial setting:

\begin{df}
The complexity of a simplicial complex $X$ with bounded degrees is the number of its simplices.
\end{df}

It is easy to see that these two notions do not change by much under transitions
$$\{\text{smooth manifolds}\}\leftrightarrows \{\text{simplicial manifolds}\}$$
given by (appropriate) triangulation, and by smoothing out the simplicial metric when possible. Therefore we don't distinguish between these two notions of complexity.

The goal of this section is to introduce some quantitative bounds on the complexity of spheres hard to cut. The original construction of~\cite{sphere} can be thought of as follows. One starts with a hard to cut graph and thickens it. It is a high-dimensional hard to cut manifold, but it has boundary and complicated topology. So one takes a tiny cone over it to simplify the topology (and then takes the boundary if one wants a closed sphere and not a ball). The theorem~\ref{celltree} provides a way to build small and locally simple cones over spaces, and that allows to keep the complexity of the resulting spheres low.

\begin{theorem}
\label{spheresep}
For any $M,\varepsilon>0$ and $n\geqslant 3$ there is a sphere $\mathcal{S}^n(M,\varepsilon)$ such that
\begin{enumerate}
\item $|\mathcal{S}|=1$,
\item $\mathcal{S}$ has complexity $C(n) \bigr(\frac{M}{\varepsilon}\bigl)^n \log(\frac{M}{\varepsilon})$,
\item $\mathcal{S}$ is $\sim(M,\varepsilon)$-hard to cut.
\end{enumerate}
\end{theorem}

In this kind of setting, if $\varepsilon$ is some fixed (not very small) number, it is commonly phrased in terms of separators.

\begin{df}
A \underline{\upshape separator} $H^{n-1}$ of a space $X^n$ is a hypersurface splitting $X^n$ into roughly equal halves. That is, $H=\partial(\Omega)$ where $\Omega\subset X$ and 
$$\vol(\Omega)>\varepsilon\vol(X)\text{ and }\vol(X\setminus \Omega)>\varepsilon \vol(X)$$
for some fixed (not very small) $\varepsilon>0$.
\end{df}

Thus, fixing such $\varepsilon$, one gets a following consequence of theorem~\ref{spheresep}:

\begin{cor}
\label{sep}
For any $M>0$ and $n\geqslant 3$ there is a sphere $\mathcal{S}^n(M)$ of volume 1 and complexity $C(n) M^n\log(M)$ so that any its separator has volume at least $M$.
\end{cor}

Note that the best complexity we can hope for here, is $\sim M^n$. Indeed, in the simplicial setting, consider a unit volume sphere of complexity $C$ (whose simplices are then of size $d\sim C^{-1/n}$). The most possible separator in such a sphere has at most $\sim C$ simplices and hence has area at most $\sim Cd^{n-1}=C^{1/n}$. Therefore $C\gtrsim M^n$, and our example gets pretty close to this rough estimate. Same analysis for general $\varepsilon$ provides an obvious bound $C\gtrsim (M/ \varepsilon)^n$ and so the theorem~\ref{spheresep} falls short of it by a log factor $\log(M/\varepsilon)$.

\medskip

The construction of the metric in~\cite{sphere} starts with embedding into $S^n$ a following graph, that is in a similar sense hard to cut.

\begin{df}
Graph $\Gamma=(V,E)$ has \underline{\upshape expander constant} $\eta(\Gamma)>0$ if for any set of vetrices $U\subset V$ there are $\eta$ times as many edges on its boundary $\partial U\subset E$:
$$|\partial U|\geq \eta \min (|U|,|V\setminus U|).$$

Define expander constant $\eta(\{\Gamma_i\})$ of a family $\{\Gamma_i\}$ of graphs $\Gamma_i$ as minimal expander constant among all individual graphs $\eta=\inf(\eta(\Gamma_i))$.
 
A family $\{\Gamma_i\}$ is \underline{\upshape an expander family} if 
$$\eta(\{\Gamma_i\})>0\text{ and }\sup(|\Gamma_i|)=\infty.$$
\end{df}

Various examples of expander families are known; for an overview see~\cite{exp}. Let us fix such a family and refer to its elements as {\itshape expander graphs}.

\subsection{Construction}

Let $(\mathcal{B},\mathcal{C},\mathcal{X})$ be the triple of $(n+1)$-ball $\mathcal{B}$, cone $\mathcal{C}$ and a (stretched) expander graph $X\simeq \mathcal{X}$, produced by theorem~\ref{celltree}. The goal is to blow up $\mathcal{X}\subset \partial \mathcal{B}$ a bit in order to make it the most significant part of $\partial \mathcal{B}$ so that it would inherit expander properties of $X$. Let $X_R$ be $X$ thickened by $R$. In more detail, instead for vertex $v$ of $X$, $X_R$ has a (triangulated at a unit scale) ball {${v_R:\cong B^{n+1}(R)}$}, and whenever vertices $v,w$ are joined in $X$ by an edge, the balls $v_R$ and $w_R$ are glued by a shared patch of their surfaces, $B^{n}(c(n)R)\hookrightarrow \partial v_R, \partial w_R$ (the patches are chosen non-overlapping for all edges incident to a given vertex). 

Note that $X_R$ is an $(n+1)$-dimensional manifold with boundary, each vertex contributing $\sim R^n$ to the boundary. Recall that the elongated graph $\mathcal{X}\subset \mathcal{B}$ has edges of length $\sim \log(N)$. So there is an embedding $\mathcal{X}\stackrel{i_X}{\hookrightarrow} \partial X_R$ for $R^n\sim \log(N)$ (such that $i_X\bigl((v,w)\bigr)\subset v_R\cup w_R$). The embedding $i_X$ is illustrated below by its image on a piece (of surface) of $X_R$.

\medskip

\includegraphics[scale=0.4]{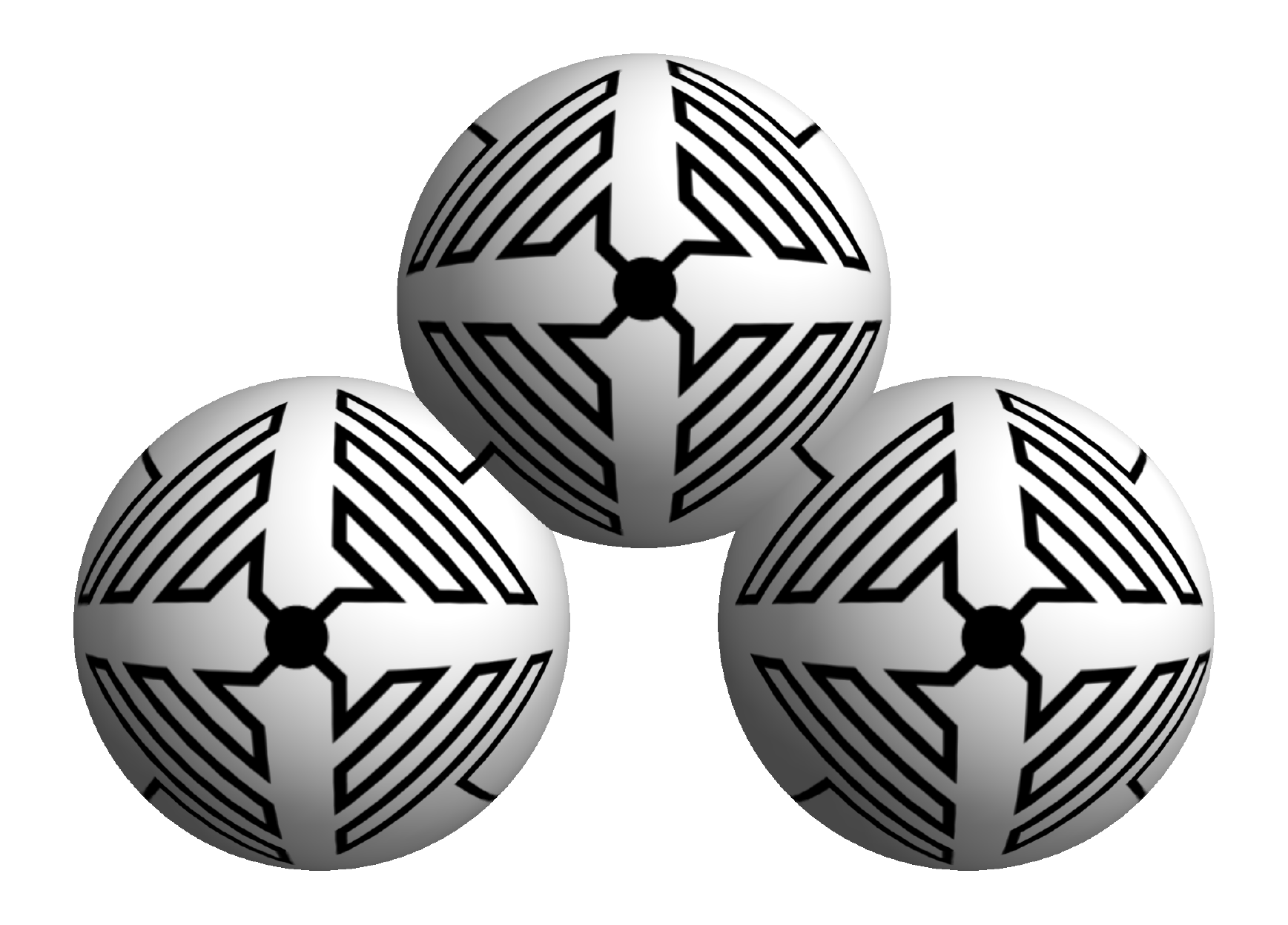}

\medskip

This embedding can be defined even on a unit tubular neighborhood of $\mathcal{X}$ in $\partial \mathcal{B}$. It can also be chosen to ``not be entangled" (isotopic to a ``straight'' embedding), so that $X_R$ retracts on $i_X(\mathcal{X})$ and therefore 
$$\mathcal{B}_R:=\mathcal{B}\cup_{i_X}X_R\simeq \mathcal{B}$$
is in a loose sense a blow-up of $\mathcal{B}$ along $\mathcal{X}$: the $\log(N)$-long and 1-thin neighborhood of $\mathcal{X}$ was replaced by $R$-thick and $R$-long version $X_R$.

Since $\mathcal{B}_R$ is an $(n+1)$-ball, its boundary is $n$-sphere $\mathcal{S}$. We show that it provides a solution to theorem~\ref{spheresep}. Specifically, given $M,\varepsilon>0$ one takes $X$ to be an expander graph on $N\sim (M/\varepsilon)^n$ vertices, then both $\mathcal{B}$ and $\partial X_R$ have $\sim N\log(N)\sim(M/\varepsilon)^n \log(M/\varepsilon)$ simplices and by construction are of bounded degrees. Hence the sphere $\mathcal{S}=\partial \mathcal{B}_R$ has these properties as well. It remains to prove that it is $\sim (M,\varepsilon)$-hard to cut when rescaled to have unit volume.

\subsection{Estimate of cut size}

Let the $\Omega\subset \mathcal{B}_R$ be the smaller part of $\mathcal{B}_R\setminus H$. By the assumption $|\Omega|>\varepsilon$, the goal is to show $|\partial \Omega|\gtrsim M$.

First assume that at least one fourth of the $\Omega$ lies in $X_R$. Then the volume of $\varepsilon$ is enough for $\sim\varepsilon(M/\varepsilon)^n$ balls of $X_R$ and by expander property of $X$, for surface 
$$\sim (\varepsilon/M)^{n-1}\dot \varepsilon(M/\varepsilon)^n=M.$$
Now if more than 3/4 of $\Omega$ are in $\mathcal{B}$, consider any simplex there. Its children in the tree of $\mathcal{B}$ are either both in $\Omega$ as well, in which case we say that they inherit both half of the original simplex (and of whatever it inherited); or at least one of them is not in $\Omega$, in which case assign this piece of $\partial \Omega$ to the considered simplex (and of whatever it inherited). Since at most 1/4 of $\Omega$ is in $X_R$ and can end inheritance chains without $\partial \Omega$, and the multiplicity of inherited material never rises above 1, at least $3/4-(1/4)-(1/4)=1/4$ 
of simplices find their way in this manner to a piece of $\partial \Omega$. Since the simplices are even smaller than the balls of $X_R$, that means that the estimate we get this way for $\partial \Omega$ is even greater then in the first case.

\begin{flushright}
$\square$
\end{flushright}


\begin{thebibliography}{20}

\bibitem{exp} O. Goldreich: Basic Facts about Expander Graphs; In: O. Goldreich, Studies in Complexity and Cryptography

\bibitem{sphere} P. Papasoglu, E. Swenson: A sphere hard to cut, arXiv:1509.02307.

\end{thebibliography}
\end{document}